\documentclass[]{interact}

\usepackage[natbibapa,nodoi]{apacite}
\theoremstyle{plain}
\newtheorem{theorem}{Theorem}[section]
\newtheorem{lemma}[theorem]{Lemma}

\theoremstyle{definition}
\newtheorem{definition}[theorem]{Definition}

\theoremstyle{remark}

\usepackage{algorithmic}
\usepackage{algorithm}
\usepackage{xcolor}

\begin{document}
  \title{Selecting Energy Efficient Inputs using Graph Structure}
  \author{
    \name{Isaac Klickstein\textsuperscript{a}\thanks{CONTACT I.~Klickstein. Email: iklick@protonmail.com} and Francesco Sorrentino\textsuperscript{a}}\affil{\textsuperscript{a}University of New Mexico, Department of Mechanical Engineering, Albuquerque NM, USA 87131}
  }
    
  \maketitle
    
  \begin{abstract}
    Selecting appropriate inputs for systems described by complex networks is an important but difficult problem that largely remains open in the field of control of networks.
    Recent work has proposed two methods for energy efficient input selection; a gradient based heuristic and a greedy approximation algorithm.
    We propose here an alternative method for input selection based on the analytic solution of the controllability Gramian of 
    {the `balloon graph', a special model graph that captures the role of both \emph{distance} and \emph{redundant paths} between a driver node and a target node.}
    The method presented is especially applicable for large networks where one is interested in controlling only a small number of outputs, or target nodes, for which current methods may not be practical because they require computing a typically very ill-conditioned matrix, called the controllability Gramian.
    Our method produces comparable results to the previous methods while being more computational efficient. 
  \end{abstract}
    
    
  \section{Introduction}
  
  Many of the systems we interact with every day are described by complex networks such as social media \citep{bovet2019influence}, the power grid \citep{arianos2009power,pagani2013power}, the world wide web \citep{barabasi2000scale}, and our own biology \citep{sporns2013structure}.
  As our ability to describe these complex networked systems improves, attention has increasingly turned to our ability to influence, or control, these systems with external signals.
  For example, targeted media campaigns (both beneficial and malicious) on social media platforms have proven to be incredibly effective \citep{grinberg2019fake}, or as our knowledge of human pharmacology grows, multi-drug multi-target therapies become viable for drug developers \citep{li2016human}.
  While the dynamics of each of these systems are significantly different, the first step toward influencing any of them requires a choice of where one should apply the external {control} signal to the system.
  {In terms of malicious social media campaigns, this means to choose which members of the social network should be targeted by counter measures to provide correct information.
  In terms of multi-drug therapies, this means to choose which drug targets to activate or inhibit by the therapy cocktail
  For power grid networks, this could mean selecting which lines should receive high voltage direct current links to improve the networks stability \cite{summers2016submodularity}.
  }
  This choice is imperative for the effectiveness of the proposed intervention.\\
  \indent
  Here, we focus on linear systems as, at least over short time scales, continuous nonlinear systems can be approximated as linear \citep{liu2016control,klickstein2017locally}.
  {Rigorous conditions for the controllability of unweighted graphs have been presented in \cite{qu2020graphical,ji2020complexity,guosufficient}. }
  The problem of selecting the smallest number of control signals to ensure a complex network is controllable has been addressed in many different frameworks such as structural controllability \citep{liu2011controllability}, exact controllability \citep{yuan2013exact}, and output controllability \citep{gao2014target,iudice2015structural,zhang2017efficient,commault2017fixed,iudice2019node}. 
  While the minimum number of inputs is sufficient to ensure controllability, applying this minimum control may lead to extremely ill-conditioned systems of equations \citep{yan2012controlling,yan2015spectrum,sun2013controllability,klickstein2017energy}.
  Instead, more recently, efficient control problems have garnered interest which look to minimize a control energy metric while constraining the number of control inputs \citep{summers2016submodularity}.
  The selection of the number of control inputs and their distribution throughout the complex network are vitally important to the feasibility and the efficiency of a control action.\\
  \indent
  Efficient controllability problems, unlike minimum controllability problems, minimize a metric on the control energy while constraining the number of control inputs \citep{summers2016submodularity}.
  Efficient control problems have previously been shown to be NP hard \citep{tzoumas2015minimal,tzoumas2016minimal} by mapping them to the hitting set problem following \citep{olshevsky2014minimal}.
  This result removes the possibility of any polynomial time algorithm to find the optimal solution.
  Instead, heuristic methods and approximation algorithms must be used to find `good', but sub-optimal, solutions.
  Two such methods are described briefly here.\\
  \indent
  The projected gradient method \citep{li2016minimum,li2018minimum,li2018enabling} finds a locally optimal solution to a continuous relaxation of the original discrete input selection problem.
  A rounding procedure, called key component analysis, is used to create a solution to the original discrete problem.
  Here, we compare our method to a version of the above heuristic which uses probabilistic projection \citep{gao2018towards} {to replace} the rounding {procedure}.\\
  \indent
  A number of control energy metrics have been shown to be submodular set functions \citep{summers2014optimal,summers2016submodularity}.
  Greedy algorithms have a well-known approximation guarantee when used to minimize submodular set functions \citep{fisher1978analysis}.
  Currently, greedy algorithms have not been explicitly used to solve the input selection problem for target control problems, except for the case that the target set coincides with the entire node set \citep{summers2016submodularity}.
  Nonetheless, the submodularity property holds for the general case of any target control problem (see corollary 2 in \citep{summers2016submodularity}) so we also compare our method to a greedy algorithm.\\
  \indent
  Both of these existing methods are iterative and require computing controllability Gramian matrices at each iteration which can be extremely ill-conditioned.
  {In this paper we present a novel method for energy efficient selection of driver nodes in general graphs. 
  As opposed to the methods described above, the method we describe here uses structural properties of the graph explicitly.
  Recent work has derived analytic expressions for the control energy of lattice networks \citep{zhao2018controllability,klickstein2018control,klickstein2018energy,klickstein2018controlb} and has shown that the control energy in these structurally simple networks can often well approximate the control energy in complex networks.
We first  analytically compute the output controllability Gramian of the `balloon graph' which consists of a number of disjoint directed paths from a driver node to a target node. This calculation is aided by the particular symmetric structure of this graph.
 Second, we solve the \emph{facility location problem} \citep{mirchandani1990discrete} with a cost matrix derived from the pair-wise costs computed using the model graph to select an energy efficient set of driver nodes.\\}
  \indent
  {In Section 2 we introduce} necessary background about graphs, the controllability Gramian, and the facility location problem.
  {In Section 3} we present our first results, deriving the controllability Gramian for a model network.
  The result is used to construct a cost matrix that describes the ability of each potential input to influence each target node.
  {In Section 4 we review two of the main alternative methods to select driver nodes from the literature. In Section 5 we present} a comparison of the three methods where we show no method out performs any other in terms of the cost of the returned solution, but that our method is more computationally efficient. {Finally, the conclusions are given in Section 6.}
  
  \section{Background}
  
  \subsection{Graph Symmetries}
  
   Graphs are denoted $\mathcal{G} = (\mathcal{V},\mathcal{E})$ which consist of $|\mathcal{V}| = n$ nodes and edges $(v_j,v_k) \in \mathcal{E}$ which may be read `from node $v_j$ to node $v_k$.'
  Unless otherwise stated, all graphs considered here are assumed to be directed.
  The set of neighbors of a node $v_j$, denoted $\mathcal{N}_j$, is defined as the set of nodes $v_k$ such that $(v_k,v_j) \in \mathcal{E}$.
  We do not include any loops, that is an edge $(v_j,v_j)$, in the set of edges, as loops are treated separately.
  A graph may be represented as an adjacency matrix, $A \in \mathbb{R}^{n \times n}$, which has elements $A_{j,k} > 0$ if $(v_k,v_j) \in \mathcal{E}$ and $A_{j,k} = 0$ otherwise.
  The diagonal of the matrix $A$, $A_{j,j} < 0$, $j = 1,\ldots,n$, represent the loops present at each node.
  In this paper, we assume \emph{uniform edge weights} and \emph{uniform loop weights}, that is, edge weights are equal, $A_{j,k} = A_{j',k'}$ for all $(v_j,v_k),(v_{j'},v_{k'}) \in \mathcal{E}$ and $(v_{j'},v_{k'}) \in \mathcal{E}$ and all loop weights are equal, $A_{j,j} = A_{k,k}$, for all $j,k = 1,\ldots,n$.
  The edge weight is denoted $\gamma > 0$ and the loop weight is denoted $-\nu < 0$.
  \begin{definition}[Graph Symmetries and the Automorphism Group \citep{lauri2016topics}]
    Let $\mathcal{G} = (\mathcal{V},\mathcal{E})$ be a graph and let $\pi : \mathcal{V} \mapsto \mathcal{V} $ be a bijection on the set of nodes of a graph.
    After applying a permutation {to the nodes of a graph}, define the new set of edges as $\mathcal{E}^{\pi}$ where if $(v_j,v_k) \in \mathcal{E}$ then $(\pi(v_j),\pi(v_k)) \in \mathcal{E}^{\pi}$.
    A permutation $\pi$ is a symmetry if $\mathcal{E} = \mathcal{E}^{\pi}$.
    The set of all such symmetries along with function composition form the automorphism group of a graph, $ \text{Aut}(\mathcal{G})$.\\
    Let $\mathcal{D} \subseteq \mathcal{V}$ be a subset of the nodes in the graph.
    The reduced automorphism group $\text{Aut}^{\mathcal{D}}(\mathcal{G})$ consists of all symmetries in $\text{Aut}(\mathcal{G})$ that do not permute any node in $\mathcal{D}$ (this concept is also known as the automorphism group of a colored graph \citep{mckay2014practical} where each driver node is a unique color and all non-driver nodes are the same color).
    \begin{equation}\label{eq:reduced_aut}
      \text{Aut}^{\mathcal{D}}(\mathcal{G}) = \left\{ \pi \in \text{Aut}(\mathcal{G}) | \pi(v_j) = v_j, \quad \forall v_j \in \mathcal{D}\right\}
    \end{equation}
  \end{definition}
  \indent
  The automorphism group (and any reduced automorphism group) induces a partition of the nodes, defined as the orbits of the graph, $\mathcal{O} = \{\mathcal{O}_1,\ldots,\mathcal{O}_q\}$, where two nodes $v_j,v_k \in \mathcal{O}_{\ell}$ if and only if there exists a symmetry $\pi$ that maps $\pi(v_j) = v_k$.
  This partition is equitable, that is, every node in orbit $\mathcal{O}_{\ell}$ has the same number of neighbors in each other orbit.
  As an example, if node $v_j$ is in orbit $\mathcal{O}_{k}$ and it has $m$ neighbors in orbit $\mathcal{O}_{\ell}$, then if node $v_{j'}$ is also in orbit $\mathcal{O}_k$ it must also have $m$ neighbors in orbit $\mathcal{O}_{\ell}$.\\
  \begin{definition}[Quotient Graph]
    Given a graph $\mathcal{G}$ and its orbits $\mathcal{O}$, the graph can be compressed to its \emph{quotient graph}, $\mathcal{Q} = (\mathcal{O},\mathcal{F})$, where each orbit is a node in the quotient graph, and the edges $(\mathcal{O}_j,\mathcal{O}_k) \in \mathcal{F}$ represent those pairs of orbits for which there exists edges passing from the nodes in $\mathcal{O}_j$ to the nodes in $\mathcal{O}_k$.
  \end{definition}
  A permutation of a set of $n$ elements, $\pi$, can also be expressed as a matrix, $P \in \{0,1\}^{n \times n}$, with elements $P_{j,k} = 1$ if $\pi(v_j) = v_k$ and $P_{j,k} = 0$ otherwise.
  Applying the permutation to the adjacency matrix yields the permuted adjacency matrix $A^{\pi} = PAP^T$.
  If $\pi$ is a symmetry then $A^{\pi} = A$ (assuming uniform edge weights and loop weights as specified above).\\
  \indent
  The orbit indicator matrix, $E \in \{0,1\}^{n \times q}$, has elements $E_{j,k} = 1$ if node $v_j \in \mathcal{O}_k$ and $E_{j,k} = 0$ otherwise.
  The adjacency matrix of the quotient graph, $A^Q \in \mathbb{R}^{q \times q}$, can be found by applying the orbit indicator matrix,
  \begin{equation}\label{eq:quotient_adj}
    A^Q = E^{\dagger} A E
  \end{equation}
  where the superscript $\dagger$ denotes the Moore-Penrose pseudoinverse, defined as $E^{\dagger} = (E^T E)^{-1} E^T$.
  The elements of the quotient graph adjacency matrix, $A_{j,k}^Q$, are equal to the number of neighbors of a node $v_{\ell} \in \mathcal{O}_j$ that are in $\mathcal{O}_k$ (multiplied by the uniform edge weight $\gamma$).
  \subsection{Minimum Energy Control}
  
  Each node is assigned a time-varying state, denoted $x_j(t)$, $j = 1,2,\ldots,n$, whose behavior is governed by its neighbors.
  We are able to influence the dynamics through a subset of the nodes, $\mathcal{D} \subseteq \mathcal{V}$, defined as the $|\mathcal{D}| = m$ \emph{driver nodes}.
  The driver node set can be represented as a matrix, $B \in \{0,1\}^{n \times m}$, where each column of $B$ has a single nonzero element corresponding to a driver node.
  An independent, external, control input, denoted $u_{\ell}(t)$, $\ell = 1,\ldots,m$, is assigned to each driver node $v_{k} \in \mathcal{D}$.
  The states evolve in time according to a system of linear differential equations where the state matrix $ A $ is the adjacency matrix of a graph.
  \begin{equation}\label{eq:dynamics}
    \dot{\textbf{x}}(t) = A \textbf{x}(t) + B \textbf{u}(t)
  \end{equation} 
  An initial condition is assigned to each node at time $t = 0$, $x_j(0) = x_{j,0}$.
  The set of $p$ target nodes, denoted $\mathcal{T} \subseteq \mathcal{V}$, are those whose states we would like to drive to a particular value at some final time $t = t_f$.
  The set of target nodes can also be represented as a matrix, $C \in \{0,1\}^{p \times n}$, where each row has a single nonzero element corresponding to a target node.
  \begin{equation}\label{eq:output}
    \textbf{y}(t) = C \textbf{x}(t)
  \end{equation}
  \begin{definition}[Controllability Gramian]
    Given matrices $A \in \mathbb{R}^{n \times n}$ and $B \in \mathbb{R}^{n \times m}$, the time-varying controllability Gramian is the symmetric $n$-by-$n$ matrix $W(t)$ that satisfies the differential Lyapunov equation (DLE),
    \begin{equation}\label{eq:dlyap}
      \dot{W}(t) = AW(t) + W(t)A^T + BB^T, \quad W(0) = O_n,
    \end{equation}
    where $O_n$ is the $n$-by-$n$ matrix of all zeroes.
    If $A$ is Hurwitz (all of its eigenvalues are in the left hand side of the complex plane) then there is a unique stable fixed point of the DLE that satisfies the algebraic Lyapunov equation (ALE),
    \begin{equation}\label{eq:lyap}
      AW + WA^T = -BB^T
    \end{equation}
    which we call the steady state controllability Gramian.
  \end{definition}
  %
  If the matrix $A$ is Hurwitz and $t_f$ is chosen large enough, then it may be appropriate to use $W$ instead of $W(t_f)$.
  \begin{lemma}[Output Controllability \citep{kailath1980linear}]\label{thm:output}
    Define the matrices $A \in \mathbb{R}^{n \times n}$, $B \in \mathbb{R}^{n \times m}$, and $C \in \mathbb{R}^{p \times n}$ and define $W(t)$ to be the solution to the DLE using $A$ and $B$.
    The triplet $(A,B,C)$ is output controllable if, for every vector $\textbf{x}_0 \in \mathbb{R}^n$, vector $\textbf{y}_f \in \mathbb{R}^p$ and positive value $t_f > 0$, there exists a time-varying signal $\textbf{u} : [0,t_f] \mapsto \mathbb{R}^m$ such that,
    \begin{equation}
      \textbf{y}_f = C e^{At_f} \textbf{x}_0 + C \int_0^{t_f} e^{A(t_f-t)} B \textbf{u}(t) dt  
    \end{equation}
    An equivalent statement is that the triplet $(A,B,C)$ is output controllable if the output controllability Gramian, 
    \begin{equation}\label{eq:output_gramian}
      \bar{W}(t) = CW(t)C^T
    \end{equation}
    is nonsingular.
  \end{lemma}
  The output controllability Gramian appears in the solution to the following optimal control problem.
  \begin{lemma}[Minimum Energy Output Control \citep{klickstein2017energy}]\label{thm:optcon}
    Define the matrices $A \in \mathbb{R}^{n \times n}$, $B \in \mathbb{R}^{n \times m}$, and $C \in \mathbb{R}^{p \times n}$ along with the vectors $\textbf{x}_0 \in \mathbb{R}^n$ and $\textbf{y}_f \in \mathbb{R}^p$.
    Then the minimum energy output control problem is,
    \begin{equation}\label{eq:optcon}
    \begin{aligned}
      \min && &J = \frac{1}{2} \int_0^{t_f} ||\bm{u}(t)||_2^2 dt\\
      \text{s.t.} && &\dot{\bm{x}}(t) = A \bm{x}(t) + B \bm{u}(t)\\
      && &\bm{x}(0) = \bm{x}_0, \quad C \bm{x}(t_f) = \bm{y}_f
    \end{aligned}
  \end{equation}
  and, if $(A,B,C)$ are output controllable, its unique solution is,
  \begin{equation}\label{eq:Jstar}
    \begin{aligned}
    J^* &= \frac{1}{2} \left( \bm{y}_f - C e^{At_f} \bm{x}_0 \right)^T \left( C W(t_f) C^T \right)^{-1} \left( \bm{y}_f - C e^{A t_f} \bm{x}_0 \right)\\
    &= \frac{1}{2} \bm{\beta}^T \bar{W}^{-1} (t_f) \bm{\beta}
    \end{aligned}
  \end{equation}
  where $\boldsymbol{\beta} = \left(\boldsymbol{y}_f - Ce^{At_f} \boldsymbol{x}_0\right)$ is called the \emph{control maneuver}, {which is the difference between the desired final output and what the final output would be in the absence of a control input}.
  {Note that the control maneuver $\bm{\beta}$ depends on the choice of the target output $\bm{y}_f$. }
  \end{lemma}
  \begin{lemma}[Symmetries in the Gramian \citep{klickstein2018control}]\label{thm:gram_symm}
  Symmetries in the graph from which the adjacency matrix $A$ was constructed appear as repeated values in the controllability Gramian.
  If two nodes $v_j,v_{j'} \in \mathcal{O}_{\ell}$ and another two nodes $v_{k}, v_{k'} \in \mathcal{O}_{\ell'}$, then the elements of the Gramian $W_{j,k}(t) = W_{j',k'}(t)$.
  \end{lemma}
  \indent
  In the following, define the driver node reduced automorphism group, $\text{Aut}^{\mathcal{D}}(\mathcal{G})$ (that is, every driver node is in an orbit of cardinality one).
  All mentions of orbit indicator matrix or quotient graph refer to those matrices and graphs induced by the driver node reduced automorphism group.\\
  \indent
  The controllability Gramian for the quotient graph satisfies Eq. \eqref{eq:dlyap} with $A$ replaced by $A^Q$ (as defined in Eq. \eqref{eq:quotient_adj}) and $B^Q = E^{\dagger} B$.
  \begin{equation}\label{eq:dlyap_quotient}
    \dot{W}^Q(t) = A^Q W^Q(t) + W^Q(t) A^{Q^T} + B^Q B^{Q^T}
  \end{equation}
  If $v_j \in \mathcal{O}_{j'}$ and $v_k \in \mathcal{O}_{k'}$ then $W_{j,k}(t) = W^Q_{j',k'}(t)$ so if $W^Q(t)$ is known, we can `expand' the result to determine $W(t)$.
  
  \subsection{Input Selection}

  Given matrices $A$, $B$, and $C$, from Lemma \ref{thm:optcon}, the minimum energy optimal control problem can be solved for $\textbf{u} : [0,t_f] \mapsto \mathbb{R}^m$ with associated cost $J^*$.
  In this paper, we are interested instead in the scenario where $A$ and $C$ are provided but we may choose the set of driver nodes $\mathcal{D} \subseteq \mathcal{V}$ (equivalently the matrix $B$ with the restrictions described previously) such that we minimize the control energy subject to a cardinality constraint on the set of driver nodes.
  For a graph with $n$ nodes, there are $\binom{n}{m}$ potential sets of $m$ driver nodes so a brute force search for even a moderate sized network is impossible.
  In addition, the solution to the minimum energy output control problem depends on the choice of control maneuver.
  To be more general, instead of minimizing the optimal cost $J^*$ in Eq. \eqref{eq:Jstar} which depends on $\bm{\beta}$ directly, an energy metric that is independent of the particular control maneuver, $M(\mathcal{D})$, is defined such that our choice of driver nodes are good in some average over the possible choices of control maneuver.
  \begin{equation}\label{eq:input_allocation}
    \begin{aligned}
      \min\limits_{\mathcal{D} \subset \mathcal{V}} && &M(\mathcal{D})\\
      \text{s.t.} && &|\mathcal{D}| = m
    \end{aligned}  
  \end{equation}
  It has previously been shown \citep{olshevsky2014minimal} that minimizing energy metrics are at least NP-hard problems, so rather than attempting to derive an algorithm to find the optimal solution to Eq. \eqref{eq:input_allocation}, heuristics and approximation algorithms must be used to return `good'  solutions (better than could be expected to be found during an extensive random search).
  The two choices of energy metric, $M(\mathcal{D})$, investigated here are the control volume and the expectation of the control energy.\\
  \indent
  The set of all control maneuvers capable of being performed with $E$ units of energy forms a $p$-dimensional ellipsoid.
  \begin{equation}\label{eq:ellipsoid}
    \mathcal{S} = \left\{\bm{\beta} \in \mathbb{R}^p | \bm{\beta}^T \bar{W}^{-1}(t_f) \bm{\beta} = E \right\}
  \end{equation}
  The volume of the ellipsoid in Eq. \eqref{eq:ellipsoid} is known to be related to the determinant of the matrix $\bar{W}(t_f)$.
  \begin{equation}\label{eq:log_volume}
    \log V(\mathcal{S}) = \log \left( \frac{\pi^{p/2}}{\Gamma (p/2+1)} \right) + \frac{1}{p} \log \det \left(\bar{W}(t_f) \right)
  \end{equation}
  The logarithm is taken of the volume as the determinant of the controllability Gramian can fall below the accuracy of double precision floating points values \citep{sun2013controllability}.
  In \citep{summers2016submodularity}, the energy metric in Eq. \eqref{eq:log_volume} is shown to be a \emph{submodular set function}.
  The submodular property of Eq. \eqref{eq:log_volume} has not been directly applied to the target control problem (only the subproblem when $p = n$).
  Nonetheless, when $p < n$, the submodularity of the log-volume holds \citep{summers2016submodularity} so we may use a greedy algorithm which retains the same approximation guarantee.
  For the metric to be minimized, we use the following inverse volume function.
  \begin{equation}\label{eq:volume_cost}
    \text{Vol}(\mathcal{D}) = -\log \det \bar{W}_{\mathcal{D}}.
  \end{equation}
  Decreasing $-\text{Vol}(\mathcal{D})$ means \emph{the set of reachable states, or control maneuvers, is larger}.\\
  %
  \indent
  The second energy metric, the expectation of the control energy over initial conditions, assumes that $\bm{y}_f = \bm{0}_p$ and $\bm{x}_0$ is a vector of $n$ independent random variables with mean zero and variance one so that $\mathbb{E}[\bm{x}_0 \bm{x}_0^T] = I_n$.
  The covariance matrix of the control maneuver can be written \citep{li2016minimum},
  \begin{equation}
    X_f = \mathbb{E}[e^{At_f} \bm{x}_0 \bm{x}_0^T e^{A^Tt_f}] = e^{At_f} e^{A^Tt_f}
  \end{equation}
  The expectation of the control energy over the control maneuvers is the metric to be minimized.
  \begin{equation}\label{eq:EB}
    \begin{aligned}
      \bar{E}(\mathcal{D}) = \text{Tr} \left(C^T \bar{W}^{-1}_{\mathcal{D}}(t_f) C X_f \right)
    \end{aligned}
  \end{equation}
  The following sections describe heuristics to solve Eq. \eqref{eq:input_allocation} when $M(\mathcal{D}) = \text{Vol}(\mathcal{D})$ and when $M(\mathcal{D}) = \bar{E}(\mathcal{D})$.
  
  \subsection{Illustrative Example}
  \begin{figure}
    \centering
    \includegraphics[scale=1]{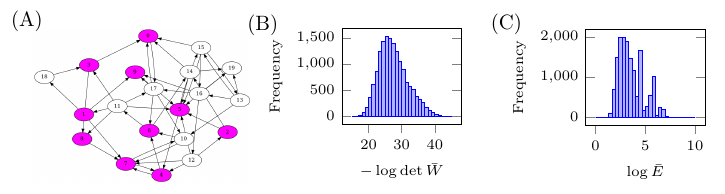}
    \caption{A small example of the input selection problem. (A) a $n = 20$ node directed graph with $p = 10 $ target nodes highlighted. (B) The volume cost from Eq. \eqref{eq:volume_cost} binned for all sets of $m = 5$ sets of nodes set as driver nodes. As the determinant of the output controllability Gramian spans over a dozen overs of magnitude, it is clear that choosing a set of driver nodes at random could lead to relatively much worse performance than the optimal solution. (C) The expectation of the control energy in Eq. \eqref{eq:EB} computed for all possible sets of $m=5$ driver nodes.}
    \label{fig:example}
  \end{figure}
  A directed graph is shown in Fig. \ref{fig:example}(A) with $n=20$ nodes and $p=10$ target nodes highlighted in pink. 
  {The goal is} to select $m=5$ driver nodes such that Eq. \eqref{eq:input_allocation} is minimized for either of the two energy metrics considered.
  This problem is small enough so that a brute force search can be employed.
  In Figs. \ref{fig:example}(B) and \ref{fig:example}(C), all sets of five nodes are set as the driver node set {successively} and both energy metrics are computed.
  In both cases, the determinant of the output Gramian and the expectation of the control energy span multiple orders of magnitude so choosing an energy efficient set of driver nodes is important.
  In the following section, we present our heuristic method to find a set of driver nodes that is energy efficient which uses only structural properties of the graph.
  \section{Graph Structure Based Input Selection}  
  
  It has previously been shown \citep{klickstein2018control,klickstein2018controlb} that the optimal cost for the single driver node and single target node problem ($m = 1$ and $p=1$) is intimately related to the structure of a graph.
  Two properties were shown to be important, the distance between the driver node and target node, $d_{j,k}$, and the number of nodes that lie along the shortest paths.
  \begin{definition}[Distance]\label{def:distance}
    A path of length $d$ is a sequence of $d$ edges, $(v_{\ell,0}, v_{\ell_1}),(v_{\ell_1},v_{\ell_2}), \ldots, (v_{\ell_{d-1}},v_{\ell_d})$.
    The distance from node $v_j$ to node $v_k$ is the length of the shortest path such $j = \ell_1$ and $k = \ell_d$.
  \end{definition}
  \begin{definition}[Redundancy]\label{def:redundancy}
    Let $\mathcal{V}_{j,k}$ be the number of nodes that lie along a shortest path from node $v_j$ to node $v_k$.
    \[
      \mathcal{V}_{j,k} = \left\{ v_{\ell} \in \mathcal{V} | d_{j,\ell} + d_{\ell,k} = d_{j,k} \right\}
    \]
    The redundancy between a pair of nodes, $v_j,v_k \in \mathcal{V}$, whose distance apart is $d_{j,k} \geq 2$, is,
    \begin{equation}\label{eq:redundancy}
      r_{j,k} = \frac{|\mathcal{V}_{j,k}| - 2}{d_{j,k} - 1}
    \end{equation}
    so that if a single path exists between two nodes, then $r_{j,k} = 1$ or if $b$ disjoint paths exist between two nodes then $r_{j,k} = b$.
  \end{definition}
  {As mentioned in the Introduction, our proposed method to select driver nodes is based on two steps. First, in Sec.\ 3.1, we analytically compute the output controllability Gramian for the directed balloon graph, from which we obtain information on the cost to control a target node at distance $d$ from a driver node with $b$ redundant paths. Then, in Sec.\ 3.2, we present a method to choose the driver nodes from the solution of the facility location problem with  a  cost  matrix derived from the pair-wise costs found for the directed balloon graph.}
  \subsection{Balloon Graph}
  \begin{figure}
    \centering
    \includegraphics[scale=1]{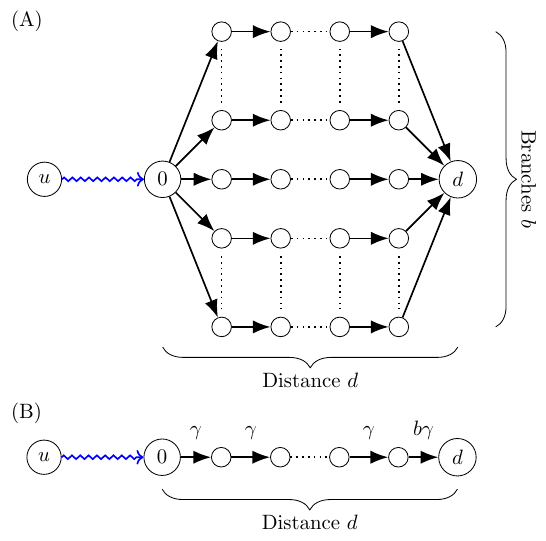}
    \caption{(A) A diagram of the directed balloon graph. There is a single driver node, labeled $v_0$, then $b$ parallel, disjoint paths of length $d$ to the terminal node, labeled $v_d$.
    All edges have uniform weight $\gamma$ and loop weight $-\nu$.
    A directed balloon graph can be completely defined by the two integers, $d$ and $b$.
    (B) A diagram of the quotient graph of the balloon graph.
    This graph is a directed path graph of uniform edge weight $\gamma$ and uniform loop weight $-\nu$, except for the right-most edge which has weight $b\gamma$.}
    \label{fig:balloon}
  \end{figure}
  A model that captures these two properties is the \emph{directed balloon graph} that consists of two end nodes, labeled $0$ and $d$, and $b$ disjoint directed paths from node $0$ to node $d$ \citep{klickstein2018energy,klickstein2018control}.
  A single control input is attached to node $0$ while $d$ is the single target node.
  Each edge is assumed to have uniform weight $\gamma > 0$ and each loop has uniform weight $-\nu < 0$.
  The driver node $0$ and target node $d$ are separated by distance $d_{0,d} = d$ and, from the definition of redundancy in Eq. \eqref{eq:redundancy}, $r_{0,d} = b$.
  The quotient graph {of the balloon graph} is a directed path graph with uniform edge weights $\gamma$ and loop weights $-\nu$ except for the last edge, $(d-1,d)$, which has edge weight $b \gamma$.\\
  \indent
  The output controllability Gramian when there is a single target node, say $\mathcal{T} = \{d\}$, is only the corresponding diagonal element $W_{d,d}(t)$.
  To determine the effect that distance and redundancy has on the control energy, this element of the controllability Gramian of the balloon graph's quotient graph is derived analytically.
  \begin{theorem}[Controllability Gramian of the Balloon Graph]\label{thm:gram_balloon}
    The diagonal element of the controllability Gramian of the balloon graph corresponding to the node $v_d$ is,
    \begin{equation}
      W_{d,d}(t) = \frac{b^2}{2\nu} \left( \frac{\gamma}{2\nu} \right)^{2d} \binom{2d}{d} \left[1 - e^{-2\nu t} \sum_{k=0}^{2d} \frac{(2\nu t)^k}{k!} \right]
    \end{equation}
    For $\nu > 0$, it can be shown that the steady state Gramian's corresponding element is,
    \begin{equation}
      W_{d,d} = \lim\limits_{t\rightarrow\infty} W_{d,d}(t) = \frac{b^2}{2\nu} \left( \frac{\gamma}{2\nu} \right)^{2d} \binom{2d}{d}
    \end{equation}
  \end{theorem}
  \begin{proof}
  Let $W_{j,k}(t)$ be the elements of the controllability Gramian of the quotient graph of the balloon graph (where the superscript $Q$ has been dropped).
  The elements $W_{j,k}(t)$ for $j,k = 0,1,\ldots,d-1$ satisfy the following system of differential equations.
    \begin{equation}\label{eq:interior}
      \begin{aligned}
        \dot{W}_{0,0}(t) &= -2\nu W_{0,0}(t) + 1\\
        \dot{W}_{j,0}(t) &= -2\nu W_{j,0}(t) + \gamma W_{j-1,0}(t), && 1 \leq j < d\\
        \dot{W}_{0,k}(t) &= -2\nu W_{0,k}(t) + \gamma W_{0,k-1}(t), && 1 \leq k < d\\
        \dot{W}_{j,k}(t) &= -2\nu W_{j,k}(t) + \gamma W_{j-1,k}(t) + \gamma W_{j,k-1}(t), && 1 \leq j,k < d
      \end{aligned}
    \end{equation}
    From symmetry, $W_{j,0}(t) = W_{0,j}(t)$, so only one set of the boundary elements in Eq. \eqref{eq:interior} must be determined.
    As every equation in Eq. \eqref{eq:interior} is linear, the Laplace transform is taken of the system where $V(s) = \mathcal{L}\{W(t)\}$.
    \begin{equation}\label{eq:balloon_sys_lyap}
      \begin{aligned}
        sV_{0,0}(s) &= -2\nu V_{0,0}(s) + \frac{1}{s}\\
        sV_{j,0}(s) &= -2\nu V_{j,0}(s) + \gamma V_{j-1,0}(s)\\
        sV_{0,k}(s) &= -2\nu V_{0,k}(s) + \gamma V_{0,k-1}(s)\\
        sV_{j,k}(s) &= -2\nu V_{j,k}(s) + \gamma V_{j-1,k}(s) + \gamma V_{j,k-1}(s)
      \end{aligned}
    \end{equation}
    The origin element, $V_{0,0}(s)$, is determined by rearranging the first line of Eq. \eqref{eq:balloon_sys_lyap}.
    \begin{equation}
      V_{0,0}(s) = \frac{1}{s(s+2\nu)}.
    \end{equation}
    The remaining elements are determined using a generating function, denoted,
    \begin{equation}
      \hat{V}(x,y;s) = \sum_{j,k = 0,1,\ldots} V_{j,k}(s) x^j y^k.
    \end{equation}
    Along the $k=0$ boundary, the elements are determined by setting $y = 0$.
    \begin{equation}
      \hat{V}(x,0;s) = \sum_{j=0,1,\ldots} V_{j,0}(s) x^j
    \end{equation}
    Multiplying the second line of Eq. \eqref{eq:balloon_sys_lyap} by $x^j$ and summing over all non-negative $j$ yields,
    \begin{equation}
      \begin{aligned}
        (s+2\nu) \sum_{j=0,1,\ldots} V_{j+1,0}(s) x^j &= \gamma \sum_{j = 0,1,\ldots} V_{j,0}(s) x^j\\
        \frac{s+2\nu}{x} \left( \hat{V}(x,0;s) - V_{0,0}(s) \right) &= \gamma \hat{V}(x,0;s)\\
        \hat{V}(x,0;s) &= \frac{s + 2\nu}{s+2\nu-\gamma x} V_{0,0}(s)
      \end{aligned}
    \end{equation}
    Define $\rho(s) = \frac{\gamma}{s+2\nu}$ so that the boundary elements can more succinctly be written as,
    \begin{equation}
      \begin{aligned}
        \hat{V}(x,0;s) &= \frac{1}{1 - \rho(s) x} V_{0,0}(s)\\
        &= V_{0,0}(s) \sum_{j \geq 0} \rho^j(s) x^j,
      \end{aligned}    
    \end{equation}
    which implies the boundary elements are,
    \begin{equation} 
      V_{j,0}(s) = \frac{1}{s(s+2\nu)} \left( \frac{\gamma}{s + 2\nu} \right)^j.
    \end{equation}
    In turn, from symmetry, the other boundary must have elements $V_{0,k}(s) = \frac{1}{s(s+2\nu)} \left( \frac{\gamma}{s+2\nu} \right)^k$.
    The interior elements are found using the two variable generating function.
    \begin{equation}
      \begin{aligned}
        \hat{V}(x,y;s) &= \frac{1}{1 - \rho(s) (x+y)} V_{0,0}(s)\\
        &= V_{0,0}(s) \sum_{\ell \geq 0} \left( \frac{\gamma}{s + 2\nu} \right)^{\ell} (x+y)^{\ell}\\
        &= V_{0,0}(s) \sum_{\ell \geq 0} \left( \frac{\gamma}{s + 2\nu} \right)^{\ell} \sum_{a = 0}^{\ell} \binom{\ell}{a} x^{\ell-a} y^a\\
        &= V_{0,0}(s) \sum_{j,k \geq 0} \binom{j+k}{k} \left( \frac{\gamma}{s+2\nu} \right)^{j+k} x^j y^k
      \end{aligned}
    \end{equation}
    The interior elements of the controllability Gramian can be read off as the $j,k$'th coefficient, $0 \leq j,k < d$,
    \begin{equation}
      V_{j,k}(s) = \frac{1}{s(s+2\nu)} \left( \frac{\gamma}{s+2\nu} \right)^{j+k} \binom{j+k}{k}
    \end{equation}
    With all elements now determined for $j,k < d$ we turn to elements when one index is equal to $d$.
    First, the element $V_{d,0}(s)$ is determined, then the elements $V_{d,j}(s)$ for $1 \leq j < d$, and then finally $V_{d,d}(s)$.
    \begin{equation}
      \begin{aligned}
        V_{d,0}(s) &= \frac{b \gamma}{s+2\nu} V_{d-1,0}(s)\\
        &= \frac{b}{s(s+2\nu)} \left( \frac{\gamma}{s+2\nu} \right)^d
      \end{aligned}
    \end{equation}
    Furthermore, it is straightforward to show that
    \begin{equation}
      V_{d,j}(s) = \frac{b}{s(s+2\nu)} \left( \frac{\gamma}{s+2\nu} \right)^{d+j} \binom{d+j}{d}
    \end{equation}
    Finally, the element of interest in the Laplace domain can be computed,
    \begin{equation}\label{eq:Vdds}
      \begin{aligned}
        V_{d,d}(s) &= \frac{b\gamma}{s+2\nu} \left(V_{d-1,d}(s) + V_{d,d-1}(s)\right)\\
        &= \frac{2b\gamma}{s+2\nu}\left(\frac{b}{s(s+2\nu)} \right) \left( \frac{\gamma}{s+2\nu} \right)^{2d-1} \binom{2d-1}{d}\\
        &= \frac{b^2}{s(s+2\nu)} \left( \frac{\gamma}{s+2\nu} \right)^{2d} \binom{2d}{d}
      \end{aligned}
    \end{equation}
    The inverse Laplace transform of $V_{d,d}(s)$ is found by using identity 5.2.18 in \citep{erdelyi1954tablesa} which states,
    \begin{equation}\label{eq:identity}
      \mathcal{L}^{-1} \left\{ \frac{1}{s(s+a)^n} \right\} = \frac{1}{a^n} \left[1 - e^{-at} \sum_{k=0}^{n-1} \frac{(at)^k}{k!}\right]
    \end{equation}
    Applying Eq. \eqref{eq:identity} to Eq. \eqref{eq:Vdds} yields the controllability Gramian element,
    \begin{equation}\label{eq:Wdd_t}
      \begin{aligned}
      W_{d,d}(t) &= \frac{b^2}{2\nu} \left( \frac{\gamma}{2\nu} \right)^{2d} \binom{2d}{d} \left[1 - e^{-2\nu t} \sum_{k=0}^{2d} \frac{(2\nu t)^k}{k!} \right]\\
      &= \frac{b^2}{2\nu} \left( \frac{\gamma}{2\nu} \right)^{2d} \frac{(2d)!}{(d!)^2} \left[1 - r(t) \right]
      \end{aligned}
    \end{equation}
    As there is a single target node, the minimum control energy in Eq. \eqref{eq:Jstar} can be written,
    \begin{equation}\label{eq:balloon_Jstar}
      J^* = \frac{\beta^2}{2} \frac{1}{W_{d,d}(t)}
    \end{equation}
    Plugging Eq. \eqref{eq:Wdd_t} into Eq. \eqref{eq:balloon_Jstar} completes the proof.
  \end{proof}
  From Eq. \eqref{eq:Jstar}, the minimum control energy for the balloon graph is $J^* \propto W_{d,d}^{-1}(t_f)$, or if $t_f$ is large enough, then $J^* \propto W_{d,d}^{-1}$.
  Using Stirling's approximation for the binomial coefficient, the control energy is approximately,
  \begin{equation}
    W_{d,d}^{-1} \approx \sqrt{\pi d} \frac{2\nu}{b^2} \left( \frac{\nu}{\gamma} \right)^d 
  \end{equation}
  Our structure based metric that approximates the control energy uses this pair-wise energy cost.
  Given a set of $m$ driver nodes $\mathcal{D}$ and a set of $p$ target nodes $\mathcal{T}$, and pairwise distances $d_{j,k}$ and redundancies $r_{j,k}$ from each node to each target node, we can construct the pair-wise cost matrix $F \in \mathbb{R}^{n \times p}$ which has elements,
  \begin{equation}\label{eq:Fjk}
    F_{j,k} = \log \sqrt{\pi d_{j,k}} \frac{2\nu}{r_{j,k}^2} \left( \frac{\nu}{\gamma} \right)^{d_{j,k}}.
  \end{equation}
  Each term, $F_{j,k}$ in Eq. \eqref{eq:Fjk} can be thought of as the cost of controlling the $k$'th target node with node $v_j$.
  The structure based metric then assigns target nodes to driver nodes by selecting which node $v_j \in \mathcal{D}$ can control the $k$'th target node the most cheaply (that is, $F_{j,k}$ is minimized over all other possible choices of driver node).
  \begin{equation}\label{eq:flp_cost}
    \text{FLP}(\mathcal{D}) = \sum_{k=1}^p \min\limits_{v_j \in \mathcal{D}} F_{j,k}
  \end{equation}
  We propose that this structure based metric can be a surrogate function to replace Eq. \eqref{eq:volume_cost} or Eq. \eqref{eq:EB} when trying to determine a set of driver nodes that minimizes one of the energy metrics.
  To test this proposal, we pick random sets of nodes from graphs using a hill climbing procedure to sample the full range of $\text{FLP}(\mathcal{D})$ and compute both $\text{FLP}(\mathcal{D})$ and $\text{Vol}(\mathcal{D})$ (or $\bar{E}(\mathcal{D})$).
  \begin{figure}
    \centering
    \includegraphics[scale=1]{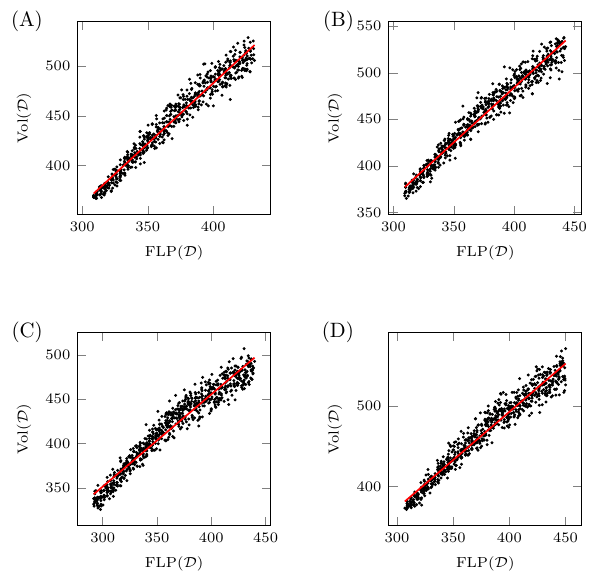}
    \caption{Comparison of the volume based cost, $\text{Vol}(\mathcal{D})$ in Eq. \eqref{eq:volume_cost}, and the structure based metric $FLP(\mathcal{D})$.
    The graphs used for the analysis are (A) a k-regular graph with $\kappa=10$, (B) an Erd\H{o}s-R\'{e}nyi graph with $\kappa_{av} = 10$, (C) a Watts-Strogatz graph with average degree $\kappa_{av} = 8$ and (D) a graph with a power law degree distribution with exponent $\gamma = 3$ and average degree $\kappa_{av} = 10$ created using the configuration model.
    All graphs are directed with $n=300$ nodes.
    The set of $p=100$ targets are chosen {uniformly at random}.
    The set of $m=33$ driver nodes are determined using the hill climbing process to achieve a desired value of $FLP(\mathcal{D})$.
    The four plots shown here are typical of all graphs examined, directed and undirected.}
    \label{fig:greedy_comp}
  \end{figure}
  \begin{figure}
    \centering
    \includegraphics[scale=1]{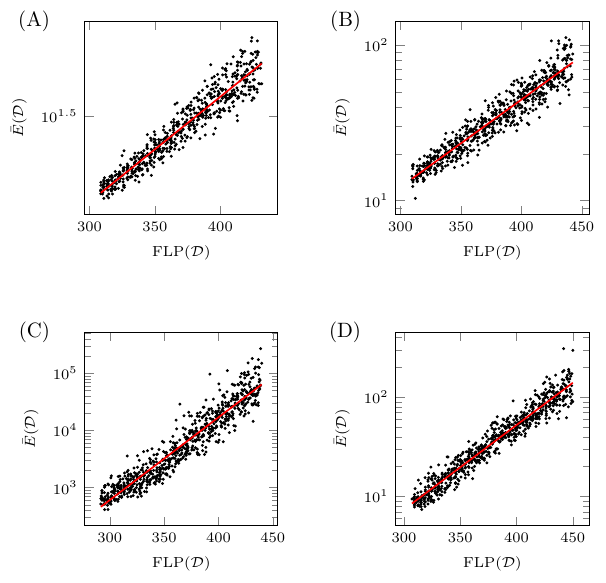}
    \caption{Comparison of the expectation energy cost in Eq. \eqref{eq:EB} and the structure based metric in Eq. \eqref{eq:flp_cost}.
    The graphs used for the analysis are (A) a k-regular graph with $\kappa=10$, (B) an Erd\H{o}s-R\'{e}nyi graph with $\kappa_{av} = 10$, (C) a Watts-Strogatz graph with average degree $\kappa_{av} = 8$ and (D) a graph with a power law degree distribution with exponent $\gamma = 3$ and average degree $\kappa_{av} = 10$ created using the configuration model.
    Each graph has $n=300$ nodes, $p= 100$ targets chosen {uniformly at random}, and $m = 33$ driver nodes chosen the same way as in Fig. \ref{fig:greedy_comp}.}
    \label{fig:lpgm_comp}
  \end{figure}
  The relationship between the ellipsoid volume cost in Eq. \eqref{eq:volume_cost} and the FLP cost in Eq. \eqref{eq:flp_cost} is shown for four example graphs in Fig. \ref{fig:greedy_comp}.
  Each graph has $n=300$ nodes, $p=100$ target nodes selected randomly, and a driver node set of $m=33$ nodes to be determined.
  The four graphs' method of construction is described in the captions of Fig. \ref{fig:greedy_comp}.
  The log volume cost in Eq. \eqref{eq:log_volume} appears on the vertical axis of each plot while the FLP cost appears on the horizontal axis.
  From the trends in Fig. \ref{fig:greedy_comp}, it is clear that if we were to find a driver node set $\mathcal{D}$ that minimized $\text{FLP}(\mathcal{D})$, that set of driver nodes would also be a competitive solution for the original optimization problem minimizing $\text{Vol}(\mathcal{D})$.\\
  \indent
  The pair-wise cost in Eq. \eqref{eq:flp_cost} is also shown to correlate with the expectation cost used by LPGM in Eq. \eqref{eq:EB}.
  A demonstration of this relation is shown in Fig. \ref{fig:lpgm_comp} for four types of graphs described in the caption.
  The two costs, $\bar{E}(\mathcal{D})$ and $\text{FLP}(\mathcal{D})$ are positively correlated as shown by the linear fitted line in red.
  Again, any driver node set $\mathcal{D}$ that minimizes $\text{FLP}(\mathcal{D})$ would be a competitive solution for the $\bar{E}(\mathcal{D})$  minimization problem as well.\\
  \indent
  Next, we present {a method based on the facility location problem \citep{mirchandani1990discrete}} to minimize $\text{FLP}(\mathcal{D})$ so that we may compare the obtained solutions with those generated by published heuristics to optimize Eq. \eqref{eq:input_allocation} with either $M(\mathcal{D}) = \text{Vol}(\mathcal{D})$ or $M(\mathcal{D}) = \bar{E}(\mathcal{D})$.
  \subsection{Facility Location Problem}
  \indent
  Facility location problems (FLP) originally arose to address the problem of choosing distribution centers to accommodate demands while minimizing transportation costs \citep{mirchandani1990discrete}.
  Let there be $p$ locations that must be supplied from $m$ distribution centers selected from $n \geq m$ possible choices.
  The cost of supplying the $j$'th location from the $k$'th distribution center is denoted $c_{j,k}$.
  Each location is assumed to be supplied from a single distribution center.
  Let the binary variables $Y_j \in \{0,1\}$, $j = 1,\ldots,n$, be the possible distribution centers where $Y_j = 1$ if it is chosen to be a distribution center and $Y_j = 0$ otherwise.
  Let the binary variables $Z_{j,k} \in \{0,1\}$, $j = 1,\ldots,n$, $k=1,\ldots,p$, denote assignments so that if distribution center $j$ supplies location $k$ then $Z_{j,k} = 1$ and $Z_{j,k} = 0$ otherwise.\\
  \indent
  The FLP can be posed as an integer linear programming (ILP) with binary variables.
  \begin{equation}\label{eq:flp}
    \begin{aligned}
      \min && &\sum_{j=1}^n \sum_{k=1}^p Y_j Z_{j,k} c_{j,k}\\
      \text{s.t.} && &\sum_{j=1}^n Y_j = m\\
      && &\sum_{j=1}^n Z_{j,k} = 1, \quad k = 1,\ldots,p\\
      && &Z_{j,k} \leq Y_j, \quad j = 1,\ldots,n, \quad k = 1,\ldots,p
    \end{aligned}
  \end{equation}
  The first constraint ensures that precisely $m$ locations are chosen to be distribution centers.
  The second constraint ensures that each location to be supplied is assigned to a single distribution center.
  The third constraint ensures locations to be supplied are only assigned to distribution centers that are opened.\\
  \indent
  Even large instances ($n \approx 1000$) of the can be solved efficiently with ILP solvers such as the GNU Linear Programming Kit \citep{glpk}.
  For larger instances {of Eq. \eqref{eq:flp}}, one can use recently developed specialized algorithms to approximately solve the FLP with an approximation guarantee \citep{jain2002new} efficiently.
  \section{Alternative Methods}
  \subsection{Greedy Algorithm}
  A greedy algorithm {that starts with an empty set and} at each iteration adds the single node to the driver node set that improves the cost function the most has an approximation guarantee of $63\%$ when the cost function is submodular \citep{fisher1978analysis}.\\
  \indent
  By the definition of the matrix $B$ we impose, the matrix product can be decomposed into the individual contributions of each driver node $BB^T = \sum_{k \in \mathcal{D}} \bm{e}_k \bm{e}_k^T$ where $\bm{e}_k$ is the unit vector with the single non-zero element corresponding to each driver node.
  This decomposition can be used to split the differential Lyapunov equation in Eq. \eqref{eq:dlyap} into the contribution of each driver node as well.
  \begin{equation}\label{eq:gram_sum}
    \begin{aligned}
      \dot{W}_k(t) &= A W_k(t) + W_k A^T + \bm{e}_k \bm{e}_k^T, \quad W_k(0) = O_n\\
      W(t) &= \sum_{k \in \mathcal{D}} W_k(t)
    \end{aligned}
  \end{equation}
  A greedy algorithm to minimize $\text{Vol}(\mathcal{D})$ over the powerset of the nodes could be applied directly assuming perfect arithmetic.\\
  \indent
  The difficulty of applying the greedy algorithm directly arises in two ways.
  The first difficulty is that storing all potential contributions of each driver node requires $np^2$ double precision variables which, if $p$ is large, could be prohibitive.
  The second difficulty is computing $\text{Vol}(\mathcal{D})$ for the first few driver node sets as the Gramian is known to have extremely small (below double precision accuracy) eigenvalues when the number of target nodes is large relative to the number of driver nodes \citep{klickstein2017energy}.
  A proposed method \citep{summers2016submodularity} to handle the first few driver nodes replaces the evaluation of $\text{Vol}(\mathcal{D})$ with $- \text{rank}_{num} (\mathcal{D})$ where the function $\text{rank}_{num}(\cdot )$ computes the \emph{numerical rank} of the output Gramian $\bar{W}_{\mathcal{D}}(t_f)$ \citep{sun2013controllability}.
  This substitute is used until enough driver nodes have been added by the greedy algorithm to ensure the output controllability Gramian is of full numerical rank.
  Algorithm \ref{alg:greedy} in Appendix \ref{apx:greedy} shows this modified version where a flag is used to perform the switch from computing the rank to the determinant.
  \subsection{$L_0$-constrained Projected Gradient Method}
  To minimize the expected energy cost in Eq. \eqref{eq:EB}, a continuous relaxation step is introduced so that the previous restrictions on $B$ are removed, that is, now $B \in \mathbb{R}^{n \times m}$.
    The main result in \citep{gao2018towards} that allows a gradient descent method to be used is the derivative of Eq. \eqref{eq:EB} with respect to $B$.
  \begin{equation}\label{eq:DEDB}
    \begin{aligned}
      \frac{\partial E(B)}{\partial B} &= -2 \int_0^{t_f} e^{A^Tt} C^T \bar{W}^{-1}_B(t_f) C X_f\\
      &\times C^T \bar{W}^{-1}_B(t_f) C e^{At} dt B
    \end{aligned}  
  \end{equation}
  With information about the gradient, a {projected gradient method} can be used {such that at each iteration the next $B$ matrix} moves in the steepest descent direction until a local minimum is found.
  A probabilistic projection is used, $\mathcal{P} : \mathbb{R}^{n \times m} \mapsto 2^{\mathcal{V}}$, that finds a set of nodes of cardinality $m$ from a dense matrix, which is the solution returned for the original optimization problem.
  Details of the algorithm can be found in Algorithm \ref{alg:projection} in Appendix \ref{apx:lpgm}.\\
  \section{Comparison}
  \subsection{Comparison with Greedy Algorithm}
  \begin{figure}
    \centering
    \includegraphics[scale=1]{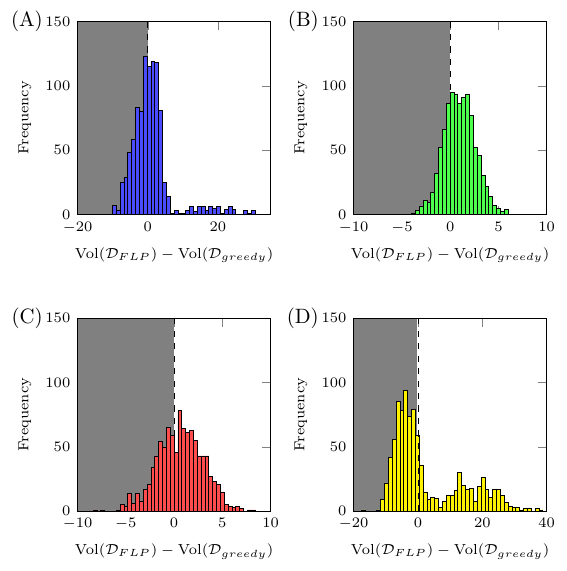}
    \caption{Comparison of the performance of the FLP formulation with the greedy algorithm.
    Each panel computes the difference between $\text{Vol}(\mathcal{D}_{FLP})$ and $\text{Vol}(\mathcal{D}_{greedy})$ defined in Eq. \eqref{eq:volume_cost}.
    The four types of graphs used are (A) Erd\H{o}s-R\'{e}nyi graphs with $\kappa_{av} = 6$, (B) $k$-regular graphs with $\kappa = 5$, (C) Watts-Strogatz graphs \citep{watts1998collective}  with $p = 5\%$, and (D) graphs with a power-law distribution with exponent $\gamma = 3$ and $\kappa_{av} = 6$.
    Each graph has $n = 50$ nodes, $p = 20$ targets and $m = 10$ driver nodes are selected.
    All graphs are undirected.
    Each panel finds the set of driver nodes returned by the FLP formulation and the greedy algorithm and compares the returned cost.
    The grey background represents cases when the FLP formulation performs better and the white background represents cases when the greedy algorithm performs better.}
    \label{fig:greedy_hist}
  \end{figure} 
  \indent
  To compare the FLP formulation described above and the greedy algorithm, we create $1000$ graphs and compute the set of driver nodes returned by the greedy algorithm in Algorithm \ref{alg:greedy} and by solving the FLP in Eq. \eqref{eq:flp}.
  In Fig. \ref{fig:greedy_hist}, $1000$ graphs of the following types are used to make the comparison; \ref{fig:greedy_hist}(A) a k-regular graph with $\kappa=10$, \ref{fig:greedy_hist}(B) an Erd\H{o}s-R\'{e}nyi graph with $\kappa_{av} = 10$, \ref{fig:greedy_hist}(C) a Watts-Strogatz graph with average degree $\kappa_{av} = 8$ and \ref{fig:greedy_hist}(D) a graph with a power law degree distribution with exponent $\gamma = 3$ and average degree $\kappa_{av} = 10$ created using the configuration model.
  Each graph is undirected and is constructed with $n = 50$ nodes and $p = 20$ nodes are chosen randomly to be in the target node set $\mathcal{T}$.
  We look for a set of $m = 10$ driver nodes, $\mathcal{D}$, such that the cost function in Eq. \eqref{eq:log_volume} is minimized.
  The set of driver nodes returned using the FLP formulation, denoted $\mathcal{D}_{FLP}$, and the set of driver nodes returned by the modified greedy algorithm, denoted $\mathcal{D}_{greedy}$, are found for each graph and their costs are computed.
  The difference of their costs,
  \begin{equation}
    D_{greedy} = \log \det (\bar{W}_{FLP}^{-1}) - \log \det (\bar{W}_{greedy}^{-1})
  \end{equation}
  is taken so that if $D_{greedy} < 0$, $\mathcal{D}_{FLP}$ is more efficient while if $D > 0$, $\mathcal{D}_{greedy}$ is more efficient.
  In Fig. \ref{fig:greedy_hist}, the cases when $\mathcal{D}_{FLP}$ is more energy efficient are shown with a gray background while the cases when $\mathcal{D}_{greedy}$ is more energy efficient are shown with a white background.
  We see that for some graph types (panels \ref{fig:greedy_hist}(A) and \ref{fig:greedy_hist}(D)), the FLP method performs better than the greedy algorithm more often, while for other graph types, the greedy algorithm performs better more often.
  Also, especially for the graphs with a power-law degree distribution in Fig. \ref{fig:greedy_hist}(D), the FLP method may not perform well as seen by the second peak in the section of the plot with a white background.\\
  \indent
  Despite the mixed results in Fig. \ref{fig:greedy_hist}, the main benefit is that our approach avoids the difficulty of computing the determinant of an ill-conditioned matrix. 
  Also, in Appendix \ref{apx:scaling}, we discuss how the greedy algorithm's computational complexity scales as $\mathcal{O}(nmp^3+n^4)$.
  To estimate the complexity of the FLP method, we use the number of nonzero {entries} in the constraint matrix, which is $(n+3np)$.
  This difference is seen in the computation times for the two methods, with the FLP solved considerably faster than the greedy algorithm.
  This means one can use both the FLP method and the greedy algorithm and take whichever solution returned has a smaller cost without increasing the amount of computational time appreciably while preserving the approximation guarantee of the greedy algorithm.

  \subsection{Comparison with LPGM}
  \begin{figure}
    \centering
    \includegraphics[scale=1]{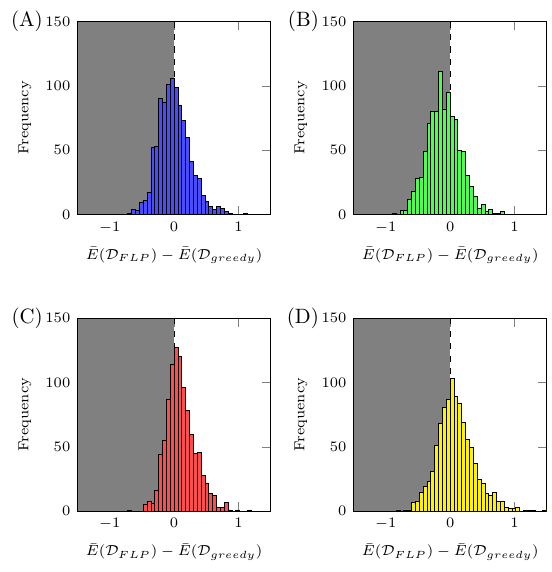}
    \caption{A comparison of the performance between LPGM and FLP method for the expectation energy cost in Eq. \eqref{eq:EB} for $1000$ realizations of the following four types of graphs.
    The four types of graphs used are (A) Erd\H{o}s-R\'{e}nyi graphs with $\kappa_{av} = 6$, (B) $k$-regular graphs with $\kappa = 5$, (C) Watts-Strogatz graphs \citep{watts1998collective}  with $p = 5\%$, and (D) graphs with a power-law distribution with exponent $\gamma = 3$ and $\kappa_{av} = 6$.
    All graphs have $n=50$ nodes, $p=20$ target nodes selected randomly, and $m=10$ driver nodes are selected.}
    \label{fig:lpgm_hist}
  \end{figure}
  A comparison of the performance of the FLP method and the LPGM heuristic for sets of four types of graphs is shown in Fig. \ref{fig:lpgm_hist}. 
  As in Fig. \ref{fig:greedy_hist}, bars in front of the gray background represent cases where the FLP algorithm returns more energy efficient driver node sets than the LPGM algorithm and vice versa for the bars with a white background.
  For the four types of graphs examined, we see that neither the FLP method nor the LPGM heuristic perform better than the other, with some slight bias towards one or the other depending on the graph.
  The benefit of the FLP method is that it scales to larger problems better than the LPGM heuristic and it does not suffer from the same overflow/underflow issues as discussed in Appendix \ref{apx:scaling}.
  \section{Conclusion}
  {The energy efficient driver node selection problem is addressed in this paper in a novel way.
  As it has previously been shown to be NP-hard, while $P\neq NP$, we cannot hope to find the optimal solution but rather we must search for `good' solutions, defined to be a solution better than that which could be reasonably expected to be found during a random search.
  While previous heuristics developed to find good solutions to this discrete optimization problem required the repeated calculation of the controllability Gramian, the method we have developed here uses the well known facility location problem with a cost matrix designed using values derived from a simple graph model.}
  {The benefits of our method are two-fold. The first is the fact that our method can provide better solutions than the previously published methods in some situations.  The second is the fact that it is efficient so it can be used in tandem with either of the previous methods without significantly increasing the computational cost.\\
  \indent
  The method presented here also exists a proof of concept that finding energy efficient sets of driver nodes can be done by using graph structure alone, rather than using properties of the controllability Gramian directly which has been shown to often be ill-conditioned or singular.
  While the cost matrix we design uses the single target single driver cost, this choice ignores the scaling of the control energy for a single driver with multiple targets.
  Future work will improve the method presented here by including terms in the cost matrix associated with a single driver assigned to multiple targets.}
  %
  %
  
  \section*{Funding}
  
  This work was supported by the National Science Foundation through grants No. 1727948 and No. CRISP-1541148.

  \section{Appendices}
  \appendix
  %
  %
    %
  %
  \section{Alternative Methods}
  Here we discuss some of the implementation details of the two alternative methods discussed in the text, namely, the greedy algorithm and the $L_0$-constrained projected gradient method.
  In both methods, to compute the controllability Gramian, we use the SLICOT routine SB03TD which is an implementation of the Bartels-Stewart algorithm.
  \subsection{Greedy Algorithm}\label{apx:greedy}
  Let $\mathcal{D}^{(k)}$ be the set of driver nodes after the $k$'th greedy step.
  The first few greedy steps correspond to the situation when only a few driver nodes have been selected so far.
  If $p$ is even of moderate size, the controllability Gramian for these first few steps will be numerically singular \citep{summers2016submodularity}, or actually singular.
  To handle this situation, the first few steps make the greedy decision based on which node increases the rank of $\bar{W}_{\mathcal{D}^{(k)}}$ the most until at some step the new controllability Gramian is of full numerical rank.
  Then the algorithm switches to choosing driver nodes corresponding to which node increases $\text{Vol}(\mathcal{D})$ the most.\\
  \indent
  To compute the rank of the matrix, a rank revealing QR factorization is performed using the SLICOT routine MB03OD.
  The determinant of a symmetric positive definite matrix is found from its Cholesky factor, $\bar{W} = LL^T$.
  Then, the determinant of $\bar{W}$ is,
  \begin{equation}
    \det \bar{W} = (\det L)^2 = \prod_{j=1}^p L_{j,j}^2
  \end{equation}
  To avoid overflow or underflow issues when taking this product, the logarithm of the determinant is computed instead.
  \begin{equation}
    \log \det \bar{W} = 2 \sum_{j=1}^p \log L_{j,j}
  \end{equation}
  We use the LAPACK routine DPOTRF to compute the Cholesky factor.
  \begin{algorithm}
    \caption{Greedy Minimization of a Set Function}\label{alg:greedy}
    \begin{algorithmic}
      \REQUIRE A desired set cardinality $m$, a state matrix $A$, a set of target nodes $\mathcal{T}$, a final time $t_f$ (possibly $\infty$).
      \FOR {$k = 1, \ldots, n$}
        \STATE Compute and store $C W_{\{k\}}(t_f) C^T$
      \ENDFOR
      \STATE Initialize $\mathcal{D} \gets \emptyset$, $\bar{W} \gets O_p$, $k \gets 0$.
      \STATE flag $\gets 1$
      \WHILE {$k < m$}
        \STATE $f_{best} \gets \infty$, $a_{best} \gets -1$
        \FOR {$j \in \mathcal{V} \backslash \mathcal{D}$}
          \IF{flag}
            \STATE $f \gets -\text{rank} (\bar{W} + CW_j(t_f)C^T)$
          \ELSE
            \STATE $f \gets -\log \det (\bar{W} + CW_j(t_f)C^T)$
          \ENDIF
          \IF {$f < f_{best}$}
            \STATE $f \gets f_{best}$, $a_{best} \gets j$
          \ENDIF
        \ENDFOR
        \STATE $\mathcal{D} \gets \mathcal{D} \cup \{a_{best}\}$
        \STATE $\bar{W} \gets \bar{W} + CW_{a_{best}}(t_f) C^T$
        \IF {flag \&\& $-C_{best} == p$}
          \STATE flag $\gets 0$
        \ENDIF
        \STATE $k \gets k+1$
      \ENDWHILE
      \RETURN $\mathcal{D}$
    \end{algorithmic}
  \end{algorithm}
  We include a flag so that the first iterations use the rank of the output controllability Gramian until at some iteration, the set of driver nodes selected so far ensures the output controllability Gramian has full \emph{numerical rank}.
  \subsection{$L_0$-constrained Projected Gradient Method}\label{apx:lpgm}
  A published algorithm proposed to solve the input selection problem to which we compare the FLP method is the $\mathcal{L}_0$-constrained projected gradient method (LPGM) \citep{gao2018towards}.
  The method combines the projected gradient method (PGM) \citep{li2016minimum,li2016optimal} which assumes all values in the $B$ matrix  with a \emph{probabilistic projection}.
  \begin{algorithm}
    \caption{Probabilistic Projection \citep{gao2018towards}}\label{alg:projection}
    \begin{algorithmic}
      \REQUIRE $B \in \mathbb{R}^{n \times m}$
      \REQUIRE $m_0 > 0$
      \STATE $r_j \gets \sum_{k=1}^m |B_{j,k}|$, $j=1,\ldots,n$
      \STATE $\mathcal{I}_{candidate} \gets \left\{j | r_j \text{ is one of the $m+m_0$ largest values} \right\}$ 
      \STATE $\mathcal{I}_{selected} \gets \emptyset$
      \WHILE {$|\mathcal{I}_{selected}| < m$}
        \STATE $p_j \gets \left\{\begin{array}{ll} r_j, & j \in \mathcal{I}_{candidate}\\ 0, & j \notin \mathcal{I}_{candidate} \end{array} \right.$, $j = 1,\ldots,n$
        \STATE Choose $j$ according to the probabilities in $p$.
        \STATE $\mathcal{I}_{selected} \gets \mathcal{I}_{selected} \cup \{j\}$
        \STATE $\mathcal{I}_{candidate} \gets \mathcal{I}_{candidate} \backslash \{j\}$
      \ENDWHILE
      \STATE $B^{L0} \gets O_{n \times m}$
      \FOR {$j \in \mathcal{I}_{selected}$ ($k = 1,\ldots,m$)}
        \STATE $B_{j,k}^{L0} \gets mr_j$
      \ENDFOR
      \STATE $B^{L0} \gets B^{L0} / \sum_{j \in \mathcal{I}_{selected}} r_j$
      \RETURN $B^{L0}$
    \end{algorithmic}
  \end{algorithm}
  The probabilistic projection in Algorithm \ref{alg:projection} appears as a step, denote $B^{L0} \gets \mathcal{P}(B)$ in the following gradient descent algorithm in Algorithm \ref{alg:gradient}.
  \begin{algorithm}
    \caption{Projected Gradient Descent \citep{gao2018towards}}\label{alg:gradient}
    \begin{algorithmic}
      \REQUIRE Graph $\mathcal{G}(\mathcal{V},\mathcal{E})$ with adj. matrix $A \in \mathbb{R}^{n \times n}$
      \REQUIRE $B_0 \in \mathbb{R}^{n \times m}$
      \REQUIRE $\mathcal{T} \subseteq \mathcal{V}$ (and corresponding $C \in \{0,1\}^{p \times n}$)
      \REQUIRE $\eta > 0$, $t_f > 0$, $K > 1$
      \STATE $E_{best} \gets \infty$
      \STATE $X_f \gets e^{A t_f} e^{A^T t_f}$
      \FOR {$k = 0,\ldots,K$}
        \STATE $B_k^{L0} \gets \mathcal{P}(B_k)$
        \STATE $W \gets \text{Lyap}(A, B_k^{L0}B_k^{L0^T}, t_f)$
        \STATE $\bar{W} \gets CWC^T$
        \STATE $E_k \gets \text{Tr} (C^T \bar{W}^{-1} C X_f)$
        \STATE $R \gets C^T \bar{W}^{-1} C X_f C^T \bar{W}^{-1} C$
        \STATE $W \gets \text{Lyap}(A^T, R, t_f)$
        \STATE $\nabla E(B_k^{L0}) \gets -2 W B_k^{L0}$
        \IF {$E < E_{best}$}
          \STATE $E_{best} \gets E$
          \STATE $B_{best}^{L0} \gets B_k^{L0}$
          \STATE $B_{k+1} \gets B_{k}^{L0}  - \eta \nabla E(B_k^{L0})$
        \ELSE
          \STATE $B_{k+1} \gets B_{k} - \eta \nabla E(B_k^{L0})$
        \ENDIF
      \ENDFOR
      \RETURN $\mathcal{D} \gets \mathcal{I}_{best}$
    \end{algorithmic}
  \end{algorithm}
  \section{Computational Cost Comparison}\label{apx:scaling}
  In the paper, namely Figs. \ref{fig:greedy_comp} and \ref{fig:lpgm_comp}, we show that the FLP cost in Eq. \eqref{eq:flp_cost} used in the ILP formulation in Eq. \eqref{eq:flp} can find competitive solutions to both the greedy algorithm with the volumetric cost in Eq. \eqref{eq:volume_cost} and the  LPGM heuristic with the expected energy cost in Eq. \eqref{eq:EB}.
  While the FLP formulation does not clearly out-perform either of the other methods in all cases, it does avoid a numerical difficulty faced by both the greedy algorithm and the LPGM heuristic.
  In the greedy algorithm, we must first compute the output controllability Gramian for each potential driver nodes' contribution, which if every node is a viable candidate, using the Bartels-Stewart algorithm \citep{bartels1972solution}, requires $\mathcal{O}(n^4)$ work.
  At each step, $k$, for $k = 1,2,\ldots,m$, we must compute either the determinant (using a Cholesky decomposition) or the rank (using a rank revealing QR decomposition) for $(n-k+1)$ $p \times p$ output controllability Gramians which both require $\mathcal{O}(p^3)$ work as we perform the comparison between each potential node to add to the set of driver nodes.
  Thus, the computational complexity of the whole greedy algorithm is $\mathcal{O}(nm p^3 + n^4)$.\\
  \indent
  The computational complexity of the LPGM heuristic, on its face, is less than the greedy algorithm, but the use of finite precision instead is the main barrier to applicability.
  To compute the descent direction, we must solve the following Lyapunov equation,
  \begin{equation}\label{eq:grad_lyap}
      O = A^T Y + Y A + R
  \end{equation}
  for the square matrix $Y$ where,
  \begin{equation}
      R = C^T \bar{W}_B^{-1} C X_f C^T \bar{W}_B^{-1} C
  \end{equation}
  can have extremely large values due to the inverse of the output controllability Gramian appearing twice.
  As Eq. \eqref{eq:grad_lyap} is a linear equation, it can also be written as $\bar{A} \cdot \text{vec}(Y) = -\text{vec}(R)$ where $\text{vec}(\cdot)$ stacks the columns of a matrix into a vector and $\bar{A} = A^T \otimes I_n + I_n \otimes A$.
  Let $|| \cdot ||$ be a vector norm, then we know that,
  \begin{equation}
      ||\bar{A} \cdot \text{vec}(Y)|| = || \text{vec}(R) || \leq || \bar{A} || \cdot || \text{vec}(Y) ||
  \end{equation}
  where $|| \bar{A}|| $ is on the order of the maximum degree in the graph so that the norm of $Y$ will be on the order of the norm of $R$, potentially very large and outside the ability of the finite precision used.
  Handling overflow issues requires care and accuracy is lost.
 The number of times this must be repeated is difficult to predict as the decay of $E_{best}$ that appears in Algorithm \ref{alg:gradient} may plateau for many iterations before decreasing \citep{gao2018towards}.\\
  \indent
  The FLP formulation as an ILP does not lend itself to an evaluation of the computational complexity directly as it depends strongly on the particular underlying algorithm and its implementation.
  An alternative metric that often correlates with the computational complexity of solving an ILP is the number of nonzero {entries} that appear in the constraint matrix.
  The constraint matrix that appears in our ILP contains $(n+3np)$ nonzeros, which grows at worst quadratically in $n$ if the number of targets $p$ grows linearly with $n$, thus it grows more slowly than the greedy algorithm.
  Our implementation which uses the GNU Linear Programming Kit \citep{glpk} to solve the ILP returns a set of driver nodes faster than the greedy algorithm every time it was compared.
  As for the LPGM, the FLP formulation does not suffer from overflow or underflow issues during the solution of the Lyapunov equation that appears in Alg. \ref{alg:gradient} to which the LPGM heuristic is prone.
  The FLP method also performed considerably faster than the LPGM heuristic for all comparisons made.
\end{document}